\documentclass[12pt]{article}
\usepackage[legalpaper,portrait,margin=1in]{geometry}
\usepackage{latexsym,amsfonts,amsmath,amssymb,mathtools,amsthm}
\usepackage{enumitem}
\setlist{listparindent=\parindent,parsep=0pt}

\usepackage{amssymb, amsmath, enumerate, hyperref,tikz,tikz-cd}
\usepackage{geometry}
\usepackage{titling}
\newcommand{\subtitle}[1]{%
  \posttitle{%
    \par\end{center}
    \begin{center}\large#1\end{center}
    \vskip0.5em}%
}
\numberwithin{equation}{section}

\usepackage[utf8]{inputenc}
\usepackage{amstext} 
\usepackage{array}   
\newcolumntype{L}{>{$}l<{$}} 

\newtheorem{theorem}{Theorem}

\newtheorem{proposition}{Proposition}
\newcommand{\RR}{\mathbb R}
\newcommand{\CC}{\mathbb C}

\newcommand{\p}{\partial}
\newcommand{\R}{\mathbb R}
\newcommand{\norm}[1]{\lVert #1 \rVert}   

\newcommand{\s}{\hspace{0.5pt}}

\newcommand{\ep}{\end{pmatrix}}

\newcommand{\ccdot}{\,\cdot\,}

\begin{document}

\title{Calder\'{o}n problem for the quasilinear conductivity equation in dimension $2$}
\author{Tony Liimatainen \and Ruirui Wu}
\date{}
\maketitle

\begin{abstract}
    In this paper we prove a uniqueness result 
    for the Calder\'{o}n problem for the quasilinear conductivity equation on a bounded domain $\R^2$. 
    The proof of the result is based on the higher order linearization method, which reduces the problem to showing  density of products of solutions to the linearized equation and their gradients. In contrast to the higher dimensional case, the proof involves delicate analysis of the correction  terms of Bukhgeim type complex geometric solutions (CGOs), which have only limited decay. To prove our results, we construct suitable families of CGOs whose phase functions have and do not have critical points. We also combine stationary phase analysis with $L^p$ estimates for the correction terms of the CGOs.
\end{abstract}

\section{Introduction}
Let $\Omega \subset \RR^2$ be a bounded open domain with smooth boundary $\p \Omega$. In this paper we consider a quasilinear conductivity equation of the form
\begin{align}\label{eq:bound_val_prob}
    \operatorname{div}(\gamma(x,u,\nabla u) \nabla u)=&0,\\
    u|_{\partial \Omega}=&f,
\end{align}
where 
\[
 \gamma : \overline{\Omega}\times \CC \times \CC^2 \rightarrow \CC
\]
is $C^\infty$ smooth function. We also assume that the quasilinear conductivity $\gamma$ satisfies
\begin{itemize}
    \item[(a)] $0<\gamma(\ccdot, 0,0) \in C^{\infty}(\bar{\Omega})$
    \item[(b)] The map $\mathbb{C} \times \mathbb{C}^2 \ni(\rho, \mu) \rightarrow \gamma(\ccdot, \rho, \mu)$ is holomorphic with values in the Hölder space $C^{1, \alpha}(\bar{\Omega})$ for some $\alpha \in(0,1)$.
\end{itemize}
With the above assumptions, the boundary value problem \eqref{eq:bound_val_prob} is well-posed in the following sense. There exists $\delta>0$ and $C>0$ such that for all 
\[
 f\in B_\delta(\partial \Omega)=\{f\in C^{2,\alpha}(\p \Omega) \ : \ \norm{f}_{C^{2,\alpha}(\p\Omega)}\}<\delta
\]
there exists a unique solution $u=u_f\in C^{2,\alpha}(\Omega)$ satisfying $\norm{u}_{C^{2,\alpha}(\overline\Omega)}<C\delta$. The proof of the above fact follows from the Banach implicit function theorem and the fact that the linearization of \eqref{eq:bound_val_prob}  is injective at the constant solution $0$. 
See for example \cite{LASSAS202144} or \cite{KKU2022partial} for similar proofs. We then define the Dirichlet-to-Neumann map $\Lambda_\gamma:B_\delta(\partial \Omega)\to C^{2,\alpha}(\p \Omega)$  by
$$
\Lambda_\gamma(f)=\left.\left(\gamma(x, u, \nabla u) \partial_\nu u\right)\right|_{\partial \Omega},
$$
where $f\in B_\delta(\partial \Omega)$, $u=u_f$, and $\nu$ is the unit outer normal to $\partial \Omega$.

We prove the following uniqueness theorem.


\begin{theorem}\label{t2}
 Let $\Omega \subset \mathbb{R}^2$ be a bounded open set with $C^{\infty}$ boundary. Assume that $\gamma_1, \gamma_2: \bar{\Omega} \times \mathbb{C} \times \mathbb{C}^2 \rightarrow \mathbb{C}$ satisfy the assumptions (a) and (b). Assume also that $\gamma_1$ and $\gamma_2$ agree up to infinite order on the boundary $\p \Omega$. Suppose that we have
$$
\Lambda_{\gamma_1}(f)=\Lambda_{\gamma_2}(f), \quad \forall f \in B_\delta(\partial \Omega).
$$
Then
$$
\gamma_1=\gamma_2 \quad \text { in } \overline{\Omega} \times \CC \times \CC^{2}.
$$
\end{theorem}
 The assumption that the unknown quantities in the theorem are known on the boundary, is to avoid proving a standard-like boundary determination results. We refer to \cite[Theorem 3.2]{carstea2024inverse} for an example of a related boundary determination result.

To prove Theorem \eqref{t2}, we use the higher order linearization method originating in the elliptic setting from \cite{LASSAS202144, FEIZMOHAMMADI20204683}. In our case the higher order linearization argument is the same as the one derived in \cite{CI} and is as follows. Denote $\gamma(\ccdot,0,0)$ by $\gamma(\ccdot,0)$. By linearization  and the uniqueness result for linear conductivity equations, we can first conclude that
\begin{align}\label{gamma_0}
   \gamma_0(\ccdot):= \gamma_1(\ccdot, 0)=\gamma_2(\ccdot, 0).
\end{align} 
Then the following can be proven by induction, see \cite{CI}¸: Assume for $k=0,1,\ldots, m-2$, that
$$
\gamma_1^{(k)}(x,0)=\gamma_2^{(k)}(x,0), \quad x\in \Omega.
$$
Then, linearizing the equation \eqref{eq:bound_val_prob} several times, one obtains
\begin{multline}
\sum_{\left(l_1, \ldots, l_m\right) \in \pi(m)} \sum_{j_1, \ldots, j_{m-1}=0}^2 \int_{\Omega}\left(\left(\partial_{\lambda_{j_1}} \ldots \partial_{\lambda_{j_{m-1}}} \gamma_1\right)(x, 0)-\left(\partial_{\lambda_{j_1}} \ldots \partial_{\lambda_{j_{m-1}}} \gamma_2\right)(x, 0)\right) \\
\left(v^{\left(l_1\right)}, \nabla v^{\left(l_1\right)}\right)_{j_1} \ldots\left(v^{\left(l_{m-1}\right)}, \nabla v^{\left(l_{m-1}\right)}\right)_{j_{m-1}} \nabla v^{\left(l_m\right)} \cdot \nabla v^{(m+1)} d x=0
\end{multline}
 for all $v^{(l)} \in C^{\infty}(\bar{\Omega})$ solving $\nabla \cdot\left(\gamma_0 \nabla v^{(l)}\right)=0$ in $\Omega, l=1, \ldots, m+1$. 

 Therefore, by also recalling the assumption (b) on $\gamma_1$ and $\gamma_2$, the proof of Theorem \ref{t2}  reduces to the following completeness result:
\begin{proposition}\label{prop:main1}
Let $\Omega \subset \mathbb{R}^2$, be a bounded open domain with $C^{\infty}$ boundary. Let $\gamma_0 \in C^{\infty}(\bar{\Omega})$ and assume that $\gamma_0$ satisfies the assumption (a). Let $m \in \mathbb{N}$ and let $T$ be $C^\infty$ smooth function with values in the space of symmetric tensors of rank $m$. Assume also that $T$ vanishes to infinite order on $\p \Omega$. Suppose that
\begin{multline}\label{prop1equ}
\sum_{\left(l_1, \ldots, l_{m+1}\right) \in \pi(m+1)} \sum_{j_1, \ldots, j_m=0}^2 \int_{\Omega} T^{j_1 \ldots j_m}(x)\left(u_{l_1}, \nabla u_{l_1}\right)_{j_1} \ldots\left(u_{l_m}, \nabla u_{l_m}\right)_{j_m} \\
\times\nabla u_{l_{m+1}} \cdot \nabla u_{m+2} \s d x=0.
\end{multline}
for all $u_l \in C^{\infty}(\bar{\Omega})$ solving $\nabla \cdot\left(\gamma_0 \nabla u_l\right)=0$ in $\Omega, l=1, \ldots, m+2$. Then $T$ vanishes identically on $\Omega$. Here $\left(u_l, \nabla u_l\right)_j, j=0,1, 2$ stands for the $j$th component of the vector $\left(u_l, \partial_{x_1} u_l, \ldots, \partial_{x_n} u_l\right)$, and $\pi(m+1)$ stands for the set of all distinct permutations of $\{1,\ldots, m +1\}$. 
\end{proposition}
To prove the above proposition, we use the reduction of the conductivity equation to the Schr\"{o}dinger equation and Bukhgeim type CGO solutions constructed for the Schr\"{o}dinger equation in \cite{Bukhgeim}. We also use modifications of the Bukhgeim solutions introduced very recently in \cite{minimal_surface_general}. These modifications include CGOs whose phase functions do not have critical points. This yields better estimates for the corresponding correction terms of the solutions, which are needed when considering inverse problems for nonlinear equations. The CGO solutions we use may have phases without critical points, or have at most one critical point. 

\subsubsection{Earlier results}
Before going into the proof of Proposition \ref{prop:main1}, we review earlier related works. For the linear conductivity equation Sylvester and Uhlmann \cite{SyU}  and Novikov \cite{Nov} proved the uniqueness result for smooth conductivities in dimensions 3 and higher. In dimension 2, Nachman \cite{Nach} proved a uniqueness result for $C^2$ conductivities. The regularity assumptions have since been relaxed by several authors. In dimensions three and higher the uniqueness is known for $W^{1, \infty}(\Omega)$ conductivities by the work of Haberman and Tataru \cite{W1inf}, and in dimension $2$ for $L^\infty$ conductivities by the work of Astala and P\"{a}iv\"{a}rinta \cite{AP}. 
 In the work by Bukhgeim \cite{Bukhgeim} the potential of a Shr\"odinger equation was recover from the corresponding DN map in dimension $2$. In the case the conductivity is matrix valued, the best results are in dimension $2$.
In the work \cite{imanuvilov2012partial} by Imanuvilov-Uhlmann-Yamamoto, a matrix valued conductivity was recover on bounded domains in $\R^2$. The very recent related work \cite{minimal_surface_general} by C\^arstea, Liimatainen and Tzou recovered the conformal structure of Riemannian surface from the Dirichlet-to-Neumann map (DN map) of the associated Shr\"odinger equation. 

The approach in the study of inverse problems for nonlinear elliptic equations was proposed in \cite{isakov1993uniqueness_parabolic}.  There the author linearized the nonlinear DN map, which reduced the inverse problem for a nonlinear equation to an inverse problem of a linear equation, which the author was able to solve by using methods for linear equations.  For the quasilinear conductivity equation, where $\gamma(x,u)$ depends on $u$, \cite{Sun} proves a uniqueness result for $C^{1,1}$ regular conductivities by linearizing of the nonlinear DN map. The linearization technique was further examined in \cite{23, 24} for elliptic equations. 
Later, second order linearizations, where data depends on two independent parameters, were used to solve inverse problems for example in  \cite{AYT2017direct,CNV2019reconstruction,KN002,sun2010inverse,sun1997inverse}.

Inverse problems for semilinear elliptic equations were also recently considered in the works by Feizmohammadi and Oksanen \cite{FEIZMOHAMMADI20204683} and Lassas, Liimatainen, Lin and Salo \cite{LASSAS202144}. These works realized how to use higher order linearizations in inverse problems for elliptic equations. The method is by now called the higher order linearization method and it was motivated by the seminal work \cite{KLU2018} by Lassas and Uhlmann that considered nonlinear hyperbolic equations. 

	After the works \cite{KLU2018, FEIZMOHAMMADI20204683,LASSAS202144}, the literature on the research of  inverse problems for nonlinear equations based on the higher order linearization method, has grown rapidly. Earlier inverse problems for quasilinear conductivity equations have been considered in \cite{hervas2002inverse, egger2014simultaneous, munoz2020calderon,shankar2020recovering}.
	The recent works \cite{LLLS2019partial,LLST2022inverse,KU2019partial,KU2019remark,feizmohammadi2023inverse, liimatainen2024uniqueness} investigated inverse problems for semilinear elliptic equations with general nonlinearities and in the case of partial data. For quasilinear conductivity equation where $\gamma(x,u, \nabla u)$, the uniqueness result we prove in this paper in dimension $2$ was obtained in \cite{CI} in dimensions $3$ and higher, and the corresponding partial data result was proven in \cite{KKU2022partial}. The works \cite{carstea2024inverse,nurminen1, nurminen2, minimal_surface_general} studied inverse problems for the minimal surface equation (which is quasilinear) on Riemannian surfaces and in Euclidean domains.

\vspace{12pt}
\textbf{Acknowledgements} R.W. would like to thank Gunther Uhlmann for proposing this project and helpful discussions throughout the progress.

T.L. was partially supported by PDE-Inverse project of the European Research Council of the
European Union, and the grant 336786 of the Research Council of Finland. Views
and opinions expressed are those of the authors only and do not necessarily reflect
those of the European Union or the other funding organizations. Neither the European Union nor the other funding organizations can be held responsible for them.

\section{CGO Solutions}\label{sec:CGOs}
In this section we construct the CGO solutions that we will use in the proof of Proposition \ref{prop:main1}.
In what follows we write $x=(x_1,x_2)\in \mathbb{R}^2 $ with $z=x_1+ix_2\in \mathbb{C}$, and $\partial = \frac{1}{2} (\partial_1-i\partial_2)$ and $\overline{\partial}= \frac{1}{2}(\partial_1+i\partial_2)$. By the identity 
\begin{equation}\label{eq1}
    -\nabla \cdot \gamma \nabla\left(\gamma^{-1 / 2} u\right)=\gamma^{1 / 2}(-\Delta+q) u,
\end{equation}
where 
\[
q=\frac{\Delta \sqrt{\gamma}}{\sqrt{\gamma}},
\] 
we obtain solutions to the conductivity equation from  solutions to the Schr\"{o}dinger equation. 

We start by recalling the construction of CGO solutions  of Bukhgeim in \cite{Bukhgeim}. For this, let $\psi$ be a holomorphic function such that $2\Im (\psi)= \varphi$ is a real valued function with one non-degenerate critical point at $z_0\in \Omega$. 
For $h>0$ and a  function $f$, we write 
$$
f_{\pm h}=e^{\pm \mathrm{i}  \varphi/h} f,
$$
and define the operators $S$ and $S^t$ by 
$$
S u=\frac{1}{4} \overline{\partial}^{-1} (b_{-h} \partial^{-1} (a_h u)),
$$
$$
S^t u = \frac{1}{4} \partial^{-1} (b_{-h} \overline{\partial}^{-1} (a_h u)),
$$
where $a\in L^p$ and $b \in W^{1,\infty}$, and 
 $\overline{\partial}^{-1}$ is the usual Cauchy operator
$$
(\overline{\partial}^{-1} u)(z)=\frac{1}{\pi} \int_{\Omega} \frac{u(w)}{z-w} \mathrm{~d} w.
$$
Here also $\partial^{-1}$ is an operator with the conjugate kernel $\pi^{-1}\overline{(z-w)^{-1}}$. 
By the proof of  in \cite[Lemma 3.1]{GuillarmouTzou}, we have 
$$
||S||_{L^r\to L^r} = O(h^{1/r}) \quad \text{and} \quad ||S||_{L^2\to L^2} = O(h^{1/2-\epsilon})
$$for any $r>2$ and $0<\epsilon<\frac{1}{2}$.

Let us then choose $a=\frac{\Delta \gamma}{\sqrt{\gamma}}$ and $b=-1$. By Theorem 3.5 in \cite{Bukhgeim} and identity \eqref{eq1}, 
\begin{align}\label{s3solu1}
    u & =\frac{1}{\sqrt{\gamma}}e^{\psi/h}(1+r_h), \\
r_h & =\sum_{n=1}^{\infty} S^n 1 \label{eq:rh}
\end{align}
is a solution to $\nabla \cdot (\gamma \nabla u)=0$. Note that the correction term $r_h$ depends on both $h$ and $\psi$. For the  convenience of the reader, we also denote
$$\bar{\partial}_{\varphi}^{-1}f := \bar{\partial}^{-1}e^{-i\varphi/h}f, \quad
\partial_{\varphi}^{-1}f := \partial^{-1}e^{i\varphi/h}f, \quad
T_h := -\partial_{\varphi}^{-1} q  \bar{\partial}_{\varphi}^{-1}
$$ 
for easier comparison to the work \cite{minimal_surface_general} from which we next borrow material from. 
With this notation, \eqref{eq:rh} reads 
$$
r_h = \bar{\partial}_{\varphi}^{-1} s_h, \quad
s_h = -\sum_{n=0}^{\infty}T_h^n \partial_{\varphi}^{-1}q.
$$
We have similarly for antiholomorphic phase function $\psi$.

In particular, after possible translation of coordinates, we can choose $ z^2$ and $-\bar{z}^2$ as the phases for the CGO solutions we use in our proofs. In these cases, the CGOs are  
\begin{align}\label{csolu1}
   {u} & =\frac{1}{\sqrt{\gamma}}e^{z^2/h}(1+{r}_h), \\
{r}_h & =\sum_{n=1}^{\infty} S^n 1
\end{align}

\begin{align}\label{csolu2}
    {u} & =\frac{1}{\sqrt{\gamma}}e^{-\overline{z}^2/h}(1+\tilde{r}_h), \\
\tilde{r}_h & =\sum_{n=1}^{\infty} (S^t)^n 1
\end{align}
By \cite[Section 4.1]{minimal_surface_general}, 
we have that the remainder  above satisfy the estimates
\begin{align}\label{cestimate}
    \left\|\tilde{r}_h\right\|_{L^r},\left\|\partial \tilde{r}_h\right\|_{L^r},\left\|\bar{\partial} \tilde{r}_h\right\|_{L^r}, ||\tilde{s}_h||_{L^r}=O\left(h^{1 / r+\epsilon_r}\right)
\end{align}
for any $r \in[2, \infty)$ and $\epsilon_r>0$ depending on $r$. We also recall from \cite[Lemma 2.2]{GuillarmouTzou} that for smooth $f$  and all 
$\epsilon>0$ small enough
\begin{equation}\label{eq:do_psi_inv}
\bar{\partial}_{\varphi}^{-1} (f) =O_{L^2}\left(h^{1/2+\epsilon}\right) \text { and } \partial_{\varphi}^{-1}(f)=O_{L^2}\left(h^{1/2+\epsilon}\right).
\end{equation}

 In addition to the above CGOs, we will also use their modification introduced recently in \cite{minimal_surface_general}.  The phase functions $\psi$ of the modifications are holomorphic or antiholomorphic, but they have no critical points in $\Omega$. 
 After possible scaling and translation of coordinates, we may assume on our domain that  $z+\frac{1}{2}z^2$ and $-\bar{z}-\frac{1}{2}\bar{z}^2$ are such phase functions. In this case we have the solutions
\begin{align}\label{ncsolu1}
   {u} & =\frac{1}{\sqrt{\gamma}}e^{(z+\frac{1}{2}z^2)/h}(1+r_h), \\
r_h & =\sum_{n=1}^{\infty} S^n 1
\end{align}
and
\begin{align}\label{ncsolu2}
   {u} & =\frac{1}{\sqrt{\gamma}}e^{(-\bar{z}-\frac{1}{2}\bar{z}^2)/h}(1+\tilde{r}_h), \\
\tilde{r}_h & =\sum_{n=1}^{\infty} (S^t)^n 1
\end{align}
to the Sh\"rodinger equation $(-\Delta+q)u=0$. 
By integration by parts we have
$$
\bar{\partial}^{-1} e^{\mathrm{i} \varphi/h} f=\frac{ih}{2}\left[e^{\mathrm{i} \varphi/h} \frac{f}{\bar{\partial} \varphi}+ \frac{ih}{2}\bar{\partial}^{-1} \left(e^{\mathrm{i}  \varphi/h} \bar{\partial}\left(\frac{f}{\bar{\partial} \varphi}\right)\right)\right],
$$
which holds for any $f \in C_0^1(\bar{\Omega})$ and  $\varphi$ having no critical points in $\Omega$. 
Therefore, using Calder\'on-Zygmund estimate (see e.g. \cite{GilbargTrudinger}), we have
$$
||\bar{\partial}^{-1} e^{\mathrm{i} \varphi/h} f||_{L^r} \leq Ch||f||_{W^{1,r}}
$$
for all $r\in (1,\infty)$. This leads to better estimates 
\begin{align}\label{ncestimate}
    ||r_h||_{L^{r}}, ||\partial r_h||_{L^{r}}, ||\bar{\partial}r_h||_{L^{r}}, ||s_h||_{L^{r}} = O(h)
\end{align}
for the correction terms in the above solutions. We refer to \cite[Section 4]{minimal_surface_general} for more details about CGOs with phases without critical points.

\section{Proof of Proposition 1}
We now prove Proposition \ref{prop:main1}, which consequently proves also Theorem \ref{t2} by the discussion in the introduction of this paper. We will choose the solutions in the integral identity \eqref{prop1equ} to be CGO solutions that we introduced in Section \ref{sec:CGOs}.  The proof is somewhat complicated and also technical. Most of the complications come from the fact that we are recovering the components of the tensor $T$ by using stationary phase, while on the other hand the correction terms of the CGOs we use satisfy $L^p$ estimates with only limited decay. It can be seen that the $L^p$ estimates are not enough to show that the correction terms  correspond to negligible terms in the asymptotic analysis of the integral identity \eqref{prop1equ}. Roughly speaking, stationary phase analysis is not well compatible with $L^p$ estimates. 

To overcome the above difficulty, we will use also CGOs whose phase functions do not have critical points, which we will see to lead to sufficiently improved decay for the correction terms. We will also use the explicit forms of the correction terms in the analysis. Since the integral identity \eqref{prop1equ} is also increasingly complicated in $m$, the proof will also be somewhat technical. For this reason, we split the proof into cases with respect to $m$. The case $m=1$ contains the most important arguments and constructions of the proof.   

\subsection{The case $m=1$:}
By assumption,  the entries $T^{j_1\cdots j_m}$ in the proposition vanish to infinite order on the boundary.
 For $m=1$, the integral identity \eqref{prop1equ} reads
\begin{align*}
0 & =\sum_{\left(l_1, l_2\right) \in \pi(2)} \sum_{j=0}^n \int T^j(x)\left(u_{l_1}, \nabla u_{l_1}\right)_j \nabla u_{l_2} \cdot \nabla u_3 \\
&=\sum_{j=0}^n \int T^j(x)\left(u_{1}, \nabla u_{1}\right)_j \nabla u_{2} \cdot \nabla u_3+\sum_{j=0}^n \int T^j(x)\left(u_{2}, \nabla u_{2}\right)_j \nabla u_{1} \cdot \nabla u_3,
\end{align*}
Where $(u_1, \nabla u_1)_j$, $j = 0, 1,\ldots, n$, stands for the $j^{\textrm{th}}$ component of  $(u_1, \p_{x_1}u_1, . . . , \p_{x_n}u_1)$, and similarly for $(u_2, \nabla u_2)_j$.
By setting $u_2=1$ in the integral identity, the terms involving $\nabla u_2$ all vanishes and so   we have
\begin{equation}\label{s3e1}
    \int  T^0(x) \nabla u \cdot \nabla v  = 0 
\end{equation}
for any functions  $u$ and $v$ solving
\begin{equation}\label{cond_equ}
    \nabla \cdot\left(\gamma_0 \nabla u\right)=0 \text { in } \Omega .
\end{equation}
We rewrite \eqref{s3e1} as
$$
\int \left(\frac{T^0}{\gamma_0}\right) \gamma_0 \nabla u \cdot \nabla v = 0
$$
We integrate by parts to move the gradient on $v$ to the other terms. By using also \eqref{cond_equ}, we get
\begin{align}\label{ibps3e1}
     \int \gamma_0 v \nabla \left(\frac{T^0}{\gamma_0}\right)\cdot \nabla u =0.
\end{align}
There is no boundary term since $T^0$ vanishes on the boundary. 
By integrating by parts again, we obtain 
\begin{align}\label{s3uvprod}
  0 =  \int \gamma_0 v \nabla \left(\frac{T^0}{\gamma_0}\right) \cdot \nabla u 
    = \int uv \nabla \cdot \left(\gamma_0 \nabla\left(\frac{T^0}{\gamma_0}\right)\right) + 
    \int \gamma_0 u \nabla \left(\frac{T^0}{\gamma_0}\right) \cdot \nabla v \nonumber \\
    = \int uv \nabla \cdot \left(\gamma_0 \nabla\left(\frac{T^0}{\gamma_0}\right)\right).
\end{align}
Here in the last identity we used \eqref{ibps3e1} with $u$ in place of $v$. 

Let us denote 
\[
A:=\nabla \cdot \left(\gamma_0 \nabla\left(\frac{T^0}{\gamma_0}\right)\right),
\]
and let $u$ and $v$ be the CGO solutions 
\begin{align*}
    u = & \frac{1}{\sqrt{\gamma_0}} e^{z^2/h} (1+r_h)\\
    v = & \frac{1}{\sqrt{\gamma_0}} e^{-\bar{z}^2/h} (1+\tilde{r}_h)
\end{align*}
as in \ref{csolu1} and \ref{csolu2}. We have by \eqref{s3uvprod} that
\begin{align}\label{s3uvsub}
\int e^{(z^2-\bar{z}^2)/h} \frac{A}{\gamma_0} (1+r_h+\tilde{r}_h+r_h\tilde{r}_h)=0.
\end{align}
By using stationary phase, we have the expansion
\[
\int e^{(z^2-\bar{z}^2)/h} \frac{A}{\gamma_0} = h\left(\frac{A}{\gamma_0}\right)(0)+o(h).
\]

(as well as for ${\partial}^{-1}$)

At this point, we mention that integration by parts also works for $\bar{\partial}^{-1}$ as it does for $\bar{\partial}$. By \cite[Theorem 1.13]{Vekua} the following holds: Let $f \in L^1$, $ \varphi \in L^{p}$ for some $p>2$, and assume that $f$ and $\varphi$ both vanish on the boundary. Then by Fubini's theorem
\begin{multline}
\int_\Omega (\bar{\partial}^{-1}f ) \varphi dz  = \frac{1}{\pi}
\int_\Omega \left( \int_{\Omega} \frac{f(\xi)}{z-\xi} d\xi \right)
 \varphi(z) dz \\
 =\frac{1}{\pi}\int_\Omega   f(\xi) \left( \int_\Omega \frac{\varphi(z)}{z-\xi}dz \right) d\xi  = -\int_\Omega f (\bar{\partial}^{-1}\varphi) dz.
\end{multline}
We use this observation for the terms in \eqref{s3uvsub}  involving the remainder terms $r_h$ and $\tilde r_h$. An integration by parts gives 
\begin{equation}\label{eq:first_integral}
\int e^{(z^2-\bar{z}^2)/h} \frac{A}{\gamma_0} r_h = 
-\int \bar{\partial}^{-1} \left(e^{(z^2-\bar{z}^2)/h} 
\frac{A}{\gamma_0}\right) e^{(\bar{z}^2-{z}^2)/h} s_h.
\end{equation}
Using \eqref{eq:do_psi_inv}  and the estimates \eqref{cestimate},
\begin{align*}
    \left\|\tilde{r}_h\right\|_{L^r},\left\|\partial \tilde{r}_h\right\|_{L^r},\left\|\bar{\partial} \tilde{r}_h\right\|_{L^r}, ||\tilde{s}_h||_{L^r}=O\left(h^{1 / r+\epsilon_r}\right), \quad r\geq 2, \quad \epsilon_r>0,
\end{align*}
we get by Cauchy's inequality that the right hand side of \eqref{eq:first_integral} is $O(h^{1+\epsilon})$ for some $\epsilon>0$. We have similarly that the terms in \eqref{s3uvsub} involving $\tilde r_h$ and $r_h\tilde r_h$  are $O(h^{1+\epsilon})$. Hence, by dividing \eqref{s3uvsub} by $h$ and letting $h\to 0$, we obtain $A(0)\gamma_0^{-1}(0)=0$. This shows that $A=0$ at $z=0$. 

By translation we can vary the critical point of the phases of the CGOs to show that $A=0$ in $\Omega$, 
$$
\nabla \cdot \left(\gamma_0 \nabla\left(\frac{T^0}{\gamma_0}\right)\right) = 0\quad \text{in $\Omega$}.
$$
Thus $T^0/\gamma_0$ is a solution to an elliptic equation in $\Omega$.  
Since $T^0$ is identically zero on the boundary, uniqueness of solutions to the Dirichlet problem of the above equation shows that 
\[
T^0=0 \text{ in } \Omega.
\]

By using that $T^0\equiv 0$, the integral identity \eqref{prop1equ} we started from now reduces to 
\begin{equation}\label{eq3}
    \sum_{\left(l_1, l_2\right) \in \pi(2)} \sum_{j=1}^2 \int T^j(x) \partial_{x_j} u_{l_1} \nabla u_{l_2} \cdot \nabla u_3 d x = 0,
\end{equation}
in the current case $m=1$ holding  for all $u_l\in {C}^{\infty}(\overline{\Omega})$, $l=1,2,3$, solving \ref{cond_equ}.
We continue by letting $u_l$, $l=1,2,3$, to  be solutions as in \eqref{ncsolu1} and \eqref{ncsolu2}: 
\begin{align*}
   {u_1} & =\frac{1}{\sqrt{\gamma_0}}e^{(z+\frac{1}{2}z^2)/h}(1+r_1), \\
   {u_2} & =\frac{1}{\sqrt{\gamma_0}}e^{(-z+\frac{1}{2}z^2)/h}(1+r_2), \\
   {u_3} & =\frac{1}{\sqrt{\gamma_0}}e^{-\bar{z}^2/h}(1+\tilde{r})
\end{align*}
Notice that we can rewrite these solutions as
\begin{align}\label{rewritesolu}
   {u_1}  =\frac{1}{\sqrt{\gamma_0}}v_1, \quad
   {u_2}  =\frac{1}{\sqrt{\gamma_0}}v_2, \quad
   {u_3}  =\frac{1}{\sqrt{\gamma_0}}v_3,
\end{align}
where the functions $v_l$ solve
$$
\Delta v +qv =0
$$
with $q = \Delta \gamma_0/\sqrt{\gamma_0}$.
Since
$$
\nabla u_{l_2} \cdot \nabla u_3 = \frac{1}{2}(\Delta (u_{l_2}u_3) - u_{l_2}\Delta u_3 - (\Delta u_{l_2})u_3),
$$
we can rewrite \eqref{eq3} as
$$
 \sum_{\left(l_1, l_2\right) \in \pi(2)} \sum_{j=1}^2 \int (T^j \partial_{x_j} u_{l_1}) (\Delta (u_{l_2}u_3) - u_{l_2}\Delta u_3 - (\Delta u_{l_2})u_3)=0
$$
By using \eqref{rewritesolu}, we also have
\begin{align}\label{Deltaid}
\Delta u_l = \Delta\left(\frac{1}{\sqrt{\gamma_0}}\right)v_l + 2\nabla\left(\frac{1}{\sqrt{\gamma_0}}\right)\cdot \nabla v_l - \frac{qv_l}{\sqrt{\gamma_0}}.
\end{align}
Let us write $\Delta=4\partial \bar{\partial}$ and integrate by parts to obtain 
\begin{align}\label{eq:second integral}
  0 &= \sum_{\left(l_1, l_2\right) \in \pi(2)} \sum_{j=1}^2 \int (T^j \partial_{x_j} u_{l_1}) (\Delta (u_{l_2}u_3) - u_{l_2}\Delta u_3 - (\Delta u_{l_2})u_3) \nonumber \\ 
  &= \sum_{\left(l_1, l_2\right) \in \pi(2)} \sum_{j=1}^2 \int \Bigg\{ 4\partial\bar{\partial}(T^j \partial_{x_j} u_{l_1}) (u_{l_2}u_3) 
  \\ 
  & \qquad \qquad-(T^j \partial_{x_j} u_{l_1})(2u_{l_2}\nabla(\frac{1}{\sqrt{\gamma_0}})\cdot \nabla v_3 + 2u_{3}\nabla(\frac{1}{\sqrt{\gamma_0}})\cdot \nabla v_{l_2}) \nonumber \\ 
   &\qquad\qquad\qquad\qquad\qquad\qquad\qquad\qquad\qquad+  (T^j \partial_{x_j} u_{l_1})F_{T,\gamma_0}(u_{l_2}v_3+v_{l_2}u_3)\Bigg\}, \nonumber
\end{align}
where $F_{T,\gamma_0}$ is a general term that depends smoothly only on $T^l$ and $\gamma_0$, and their derivatives.
As before, by estimates \eqref{ncestimate}, \eqref{cestimate} and stationary phase,
we have
\begin{equation}\label{eq:Tgamma_id}
\int (T^j \partial_{x_j} u_{l_1}F_{T,\gamma_0})(u_{l_2}v_3+v_{l_2}u_3)=O(1).
\end{equation}
Using the above and by writing   $\nabla u \cdot \nabla v = 2(\partial u \bar{\partial}v+ \bar{\partial} u \partial v)$, the integral in  \eqref{eq:second integral} reads   
\begin{align*}
&\int 4\partial\bar{\partial}(T^j \partial_{x_j} u_{l_1}) (u_{l_2}u_3) - \int 4(T^j \partial_{x_j} u_{l_1}) \left[u_{l_2}\left(\partial \left(\frac{1}{\sqrt{\gamma_0}}\right)  \bar{\partial}v_3+\bar{\partial} \left(\frac{1}{\sqrt{\gamma_0}}\right)  {\partial}v_3\right)\right.\\
&\qquad\quad\qquad\qquad\qquad\qquad\qquad\qquad \left.+u_3\left(\partial \left(\frac{1}{\sqrt{\gamma_0}}\right)  \bar{\partial}v_{l_2} +\bar{\partial} \left(\frac{1}{\sqrt{\gamma_0}}\right) {\partial}v_{l_2}\right)\right] +O(1).
\end{align*}
Here we have 
\begin{multline}
    \int 4\partial\bar{\partial}(T^j \partial_{x_j} u_{l_1}) (u_{l_2}u_3)=\int 4(\partial T^j \bar{\partial} \partial_{x_j} u_{l_1} + 4 \bar{\partial} T^j {\partial} \partial_{x_j} u_{l_1}) (u_{l_2}u_3) \\
+\int4T^j\partial_{x_j}\left(\partial \left(\frac{1}{\sqrt{\gamma_0}}\right) \bar{\partial}v_{l_1}+ \bar{ \partial} \left(\frac{1}{\sqrt{\gamma_0}}\right) {\partial}v_{l_1}\right) (u_{l_2}u_3) +O(1)
\end{multline}
since in the case both 
$\partial$ and $\bar{\partial}$ hit $T^j$ in the term $4\partial\bar{\partial}(T^j \partial_{x_j} u_{l_1})$,  then the resulting integral term will be $O(1)$ by the argument after \eqref{eq:first_integral}. We also used \eqref{Deltaid} and \eqref{eq:Tgamma_id} again.

Note next that
\[
\int \partial T^j \bar{\partial} \partial_{x_j} u_{l_1}=O(1),
\]
because the phase function of $u_{l_1}$ is holomorphic. Indeed, when $\bar{\partial}$ hits the exponential factor $\exp{((\pm z+\frac{1}{2}z^2)/h)}$ of $u_{l_1}$ the result vanishes and the situation is then similar to the case where there are only first order derivatives of $u_{l_1}$. Regarding this, we note that when $\bar\partial\partial$ hits the correction term, the corresponding integral is also $O(1)$. This is due to the Calder\'on-Zygmund estimate explained below.

By the Calder\'on-Zygmund inequality (see for example \cite[Corollary 9.10]{GilbargTrudinger}), the $L^2$ norm of any second order derivative of $r_{l_1}$ multiplied by a compactly supported $C^\infty$ function, say $H$, can be estimated as 

\begin{multline}
||H\nabla^2  r_{l_1}||_{L^2} \leq ||\nabla^2 (H r_{l_1})||_{L^2}+2 ||\nabla H\otimes \nabla r_{l_1}||_{L^2} +||r_{l_1}\nabla^2 H||_{L^2} \\
\lesssim ||\partial\bar{\partial}(Hr_{l_1})||_{L^2} + O(h^{1/2+\epsilon}) 
\lesssim 
||\partial(e^{i\varphi/h}s_h)||_{L^2} +O(h^{1/2+\epsilon}) \\
\lesssim
\left(\frac{1}{h} ||s_h||_{L^2}+||\partial s_h||_{L^2}\right)+O(h^{1/2+\epsilon})=O(1).
\end{multline}
Here we used the estimate $||s_h||_{L^2}=O(h)$ from \ref{ncestimate} and  
\[
s_h = -\sum_{k=0}^{\infty}T_h^{\s k} \partial_{\varphi}^{-1}q
\]
to have that $||\partial s_h||_{L^2}=O(1)$. So the third identity holds.


Returning to the main line of the proof, combining our computations so far shows that  \eqref{eq:second integral} equals
\begin{align}\label{s3quantity1}
&\sum_{\left(l_1, l_2\right) \in \pi(2)} \sum_{j=1}^2 \int  \Bigg\{4 \bar{\partial} T^j ({\partial} \partial_{x_j} u_{l_1}) u_{l_2}u_3\nonumber\\ 
&+4T^j\partial_{x_j}\left(\left(\frac{1}{\sqrt{\gamma_0}}\right) \bar{\partial}v_{l_1}+ \bar{ \partial} \left(\frac{1}{\sqrt{\gamma_0}}\right) {\partial}v_{l_1}\right) u_{l_2}u_3 \nonumber \\ 
&- 4(T^j \partial_{x_j} u_{l_1}) \left[u_{l_2}\left(\partial \left(\frac{1}{\sqrt{\gamma_0}}\right)  \bar{\partial}v_3+\bar{\partial} \left(\frac{1}{\sqrt{\gamma_0}}\right)  {\partial}v_3\right)\right.\\
&\qquad\quad\qquad\qquad\qquad\qquad\qquad\qquad \left.+u_3\left(\partial \left(\frac{1}{\sqrt{\gamma_0}}\right)  \bar{\partial}v_{l_2} +\bar{\partial} \left(\frac{1}{\sqrt{\gamma_0}}\right) {\partial}v_{l_2}\right)\right]\Bigg\}+O(1).\nonumber
\end{align}
Let us then recall the notations  $\partial_{x_1}=\partial_1 = \partial+\bar{\partial}$ and $ \partial_{x_2}=\partial_2= i(\partial-\bar{\partial})$. So we have
\begin{equation}\label{pp12}
\partial \partial_{1} = \partial^2+\partial \bar{\partial}, \quad
\partial \partial_{2} = i(\partial^2-\partial \bar{\partial}).
\end{equation}
Let us then consider the first term in \eqref{s3quantity1}. We have
\[
\sum_{\left(l_1, l_2\right) \in \pi(2)} \sum_{j=1}^2 \int  4 \bar{\partial} T^j ({\partial} \partial_{x_j} u_{l_1}) u_{l_2}u_3=\sum_{\left(l_1, l_2\right) \in \pi(2)}  \int  4 \bar{\partial} (T^1+iT^2) (\partial^2 u_{l_1}) (u_{l_2}u_3)+O(1)
\]
by arguing as above. 
Consequently 
\begin{align}\label{first_int_result}
&\sum_{\left(l_1, l_2\right) \in \pi(2)} \sum_{j=1}^2 \int  4 \bar{\partial} T^j ({\partial} \partial_{x_j} u_{l_1}) u_{l_2}u_3=  \sum_{\left(l_1, l_2\right) \in \pi(2)} \int  4 \bar{\partial} (T^1+iT^2) ({\partial}^2 u_{l_1}) u_{l_2}u_3 +O(1) \nonumber \\
& = \int \frac{4}{h^2 \gamma_0^{3/2}}\bar{\partial} (T^1+iT^2)
e^{(z^2-\bar{z}^2)/h}\left((1+z)^2+(-1+z)^2\right)(1+r_1)(1+r_2)(1+r_3)+O(1) \nonumber \\&
= \frac{8}{h} (\gamma_0(0))^{-3/2}\s  \bar{\partial} (T^1+iT^2) (0)+o(h^{-1})
\end{align}
by the stationary phase. Here we also used that the integral corresponding to $\partial^2 r_{l_1}$ is $O(1)$ by the Calder\'on-Zygmund estimate.  

For the second term of \eqref{s3quantity1},
 we similarly have 
\begin{align}\label{second_integral_result}
  & \sum_{\left(l_1, l_2\right) \in \pi(2)} \sum_{j=1}^2 \int  4T^j\partial_{x_j}\left[\partial \left(\frac{1}{\sqrt{\gamma_0}}\right) \bar{\partial}v_{l_1}+ \bar{ \partial} \left(\frac{1}{\sqrt{\gamma_0}}\right) \partial v_{l_1}\right] u_{l_2}u_3 \nonumber \\ 
   & =  \sum_{\left(l_1, l_2\right) \in \pi(2)}
   \int 4(T^1+i T^2) \bar{\partial}\left(\frac{1}{\sqrt{\gamma_0}}\right) \left( \partial^2 v_{l_1}\right) u_{l_2}u_3+O(1) \nonumber \\
   & = \int \frac{4}{h^2\gamma_0}(T^1+iT^2) \bar{\partial}\left(\frac{1}{\sqrt{\gamma_0}}\right) \left[ e^{(z^2-\bar{z}^2)/h} ((1+z)^2+(-1+z)^2) (1+r_1)(1+r_2)(1+r_3)\right] \nonumber \\
   &+O(1)
   = \frac{8}{h} \gamma_0(0)^{-1}\bar{\partial}(\gamma_0(0)^{-1/2}) \s   (T^1+iT^2) (0)+o(h^{-1}).
\end{align}

For the last term of \eqref{s3quantity1}, we first notice that 
\begin{align*}
   & \sum_{\left(l_1, l_2\right) \in \pi(2)} \sum_{j=1}^2 \int - 4T^j \partial_{x_j} u_{l_1} \Bigg\{ u_{l_2} \left[ \partial \left(\frac{1}{\sqrt{\gamma_0}}\right) \bar{\partial}v_3 +\bar{\partial}\left(\frac{1}{\sqrt{\gamma_0}}\right) {\partial}v_3\right] \\
&\qquad\quad\qquad\qquad\qquad\quad\qquad\qquad+u_3\left[ \partial \left(\frac{1}{\sqrt{\gamma_0}}\right)\bar{\partial}v_{l_2} +\bar{\partial} \left(\frac{1}{\sqrt{\gamma_0}}\right){\partial}v_{l_2}\right]\Bigg\}\\&
    =\sum_{\left(l_1, l_2\right) \in \pi(2)} \sum_{j=1}^2 \int - 4T^j \partial_{x_j} u_{l_1}\left[u_{l_2}  \partial \left(\frac{1}{\sqrt{\gamma_0}}\right) \bar{\partial}v_3+ u_3 \bar{\partial}\left( \frac{1}{\sqrt{\gamma_0}} \right){\partial}v_{l_2}\right] +o(h^{-1})
\end{align*}
by arguing similarly as before. Then we compute
\begin{align*}
   & \sum_{\left(l_1, l_2\right) \in \pi(2)} \sum_{j=1}^2 \int - 4(T^j \partial_{x_j} u_{l_1})u_{l_2} \partial \left(\frac{1}{\sqrt{\gamma_0}}\right) \bar{\partial}v_3 \\
    & = \int -\frac{4(T^1+iT^2)}{h^2 \gamma_0} \partial \left(\frac{1}{\sqrt{\gamma_0}}\right) \left[e^{(z^2+\bar{z}^2)/h}((1+z)  +(-1+z))(-2\bar{z})(1+r_1)(1+r_2)(1+r_3)\right] \\
    &+o(h^{-1}) = o(h^{-1}). 
\end{align*}
We also obtain 
\begin{align}\label{third_integral_result}
     &\sum_{\left(l_1, l_2\right) \in \pi(2)} \sum_{j=1}^2 \int - 4(T^j \partial_{x_j} u_{l_1}) u_3 \bar{\partial} \left(\frac{1}{\sqrt{\gamma_0}}\right) {\partial}v_{l_2} \nonumber \\& =
     \sum_{\left(l_1, l_2\right) \in \pi(2)} \int_\Omega -4(T^1+iT^2) \bar{\partial} \left(\frac{1}{\sqrt{\gamma_0}}\right) (\partial u_{l_1})(\partial v_{l_2}) u_3 + O(1) \nonumber \\ &
     = \int_\Omega \frac{-4}{h^2 \gamma_0} (T^1+iT^2) \bar{\partial} \left(\frac{1}{\sqrt{\gamma_0}}\right)\left[  e^{(z^2+\bar{z}^2)/h} 2(1+z)(-1+z)(1+r_1)(1+r_2)(1+r_3)\right] \nonumber \\
     &+ O(1) 
     =\frac{8}{h}\gamma_0(0)^{-1}\bar{\partial}(\gamma_0(0)^{-1/2}) \s  \bar{\partial} (T^1+iT^2) (0)+o(h^{-1}).
\end{align}

Combining the above the results in \eqref{first_int_result}, \eqref{second_integral_result} and \eqref{third_integral_result},  multiply the right hand side of \eqref{eq:second integral} by $h$ and letting $h\to \infty$ shows that 
$$
{\bar{\partial}(T^1+iT^2)}- \frac{\bar{\partial}{\gamma_0}}{{\gamma_0}}(T^1+iT^2)=0
$$
in $\Omega$. By applying $\partial$ to both sides of the above equation, we get a second order elliptic equation for $T^1+iT^2$. Since $T^1$ and $T^2$ are real and vanish to first order on $\partial \Omega$ by boundary determination, we conclude that $T^1+iT^2=0$ in $\Omega$ by unique continuation.

So far we have shown that $T^0=T^1+iT^2=0$. Let us next choose 
\begin{align*}
   {u_1} & =\frac{1}{\sqrt{\gamma_0}}e^{(-\bar{z}-\frac{1}{2}\bar{z}^2)/h}(1+r_1), \\
   {u_2} & =\frac{1}{\sqrt{\gamma_0}}e^{(\bar{z}-\frac{1}{2}\bar{z}^2)/h}(1+r_2), \\
   {u_3} & =\frac{1}{\sqrt{\gamma_0}}e^{{z}^2/h}(1+\tilde{r})
\end{align*}
as solutions to $\nabla \cdot\left(\gamma_0 \nabla u\right)=0$. 
By using these solutions, and arguing in a similar manner as before, 
we obtain $T^1-iT^2=0$ in $\Omega$. Combining everything, we thus have shown $T^0=T^1=T^2=0$ in $\Omega$. This concludes the proof of the case $m=1$. 
\vspace{12pt}

\subsection{The case $m=2$:}
For $m=2$, the integral identity \eqref{prop1equ} reads
\begin{equation}\label{eq:m2identity}
\sum_{\left(l_1, l_2,l_3\right) \in \pi(3)} \sum_{j_1, j_2=0}^2 \int T^{j_1j_2}(x)\left(u_{l_1}, \nabla u_{l_1}\right)_{j_1} \left(u_{l_2}, \nabla u_{l_2}\right)_{j_2} 
\nabla u_{l_{3}} \cdot \nabla u_{4}\s =0.
\end{equation}
If we let two of the functions $u_1, u_2,u_{3}$  be constant functions equal to $1$, then we get
\begin{equation*}
    \int  T^{00}(x) \nabla u \cdot \nabla v  = 0 ,
\end{equation*}
which is the same as the identity \eqref{s3e1} in  the $m=1$ case for $T^{0}$. This proves $T^{0 0}=0$ in $\Omega$. Next, we let one of $u_1, u_2,u_{3}$ be the constant function $1$. This yields
\begin{equation*}
    \sum_{\left(l_1, l_2\right) \in \pi(2)} \sum_{j=1}^2 \int T^{0j}(x) \partial_{x_j} u_{l_1} \nabla u_{l_2} \cdot \nabla u_3  = 0,
\end{equation*}
which is the same identity \eqref{eq3} as in the
 $m=1$ case for $T^1$ and $T^2$. Thus we obtain  $T^{01}$ and  $T^{02}=0$ in $\Omega$.

By using $T^{00}=T^{01}=T^{02}=0$,  the identity \eqref{eq:m2identity} becomes
\begin{align}\label{eqT^ij_neq0}
    \sum_{\left(l_1, l_2, l_3\right) \in \pi(3)} \sum_{j, k=1}^2 \int T^{j k}(x) \partial_{x_j} u_{l_1} \partial_{x_k} u_{l_2} \nabla u_{l_3} \cdot \nabla u_4 =0.
\end{align}
Next we choose solutions $u_l$, $l=1,\ldots, 4$, as in \eqref{ncsolu1}
and \eqref{csolu2}, as we did when proving $T^1+iT^2=0$ in $m=1$ case. That is, we choose
\begin{align*}
   {u_1} & =\frac{1}{\sqrt{\gamma_0}}e^{(z+\frac{1}{3}z^2)/h}(1+r_1), \\
   {u_2} & =\frac{1}{\sqrt{\gamma_0}}e^{(z+\frac{1}{3}z^2)/h}(1+r_2), \\
   {u_3} & =\frac{1}{\sqrt{\gamma_0}}e^{(-2{z}+\frac{1}{3}{z}^2)/h}(1+r_3),\\
   {u_4} & =\frac{1}{\sqrt{\gamma_0}}e^{-\bar{z}^{2}/h}(1+r_4).
\end{align*}
Rewrite \eqref{eqT^ij_neq0} as
$$
 \sum_{\left(l_1, l_2\right) \in \pi(2)} \sum_{j,k=1}^2 \int (T^{jk} \partial_{x_j} u_{l_1} \partial_{x_k} u_{l_2}) (\Delta (u_{l_3}u_4) - u_{l_3}\Delta u_4 - (\Delta u_{l_3})u_4)=0
$$
and we can write solutions as
\begin{align*}
   {u_l}  =\frac{1}{\sqrt{\gamma_0}}v_l, \qquad l=1,2,3,4.
\end{align*}
We obtain by substituting the solutions
\begin{align}
  0 &= \sum_{\left(l_1, l_2, l_3\right) \in \pi(3)} \sum_{j,k=1}^2 \int (T^{jk} \partial_{x_j} u_{l_1} \partial_{x_k}u_{l_2} )(\Delta (u_{l_3}u_4) - u_{l_3}\Delta u_4 - (\Delta u_{l_3})u_4) \nonumber \\ 
  &= \sum_{\left(l_1, l_2, l_3\right) \in \pi(3)} \sum_{j=1}^2 \int \Bigg\{ 4\partial\bar{\partial}(T^{jk} \partial_{x_j} u_{l_1}\partial_{x_k}u_{l_2}) (u_{l_3}u_4) 
  \\ 
  & \qquad \qquad-(T^{jk} \partial_{x_j} u_{l_1}\partial_{x_k}u_{l_2})(2u_{l_3}\nabla(\frac{1}{\sqrt{\gamma_0}})\cdot \nabla v_4 + 2u_{4}\nabla(\frac{1}{\sqrt{\gamma_0}})\cdot \nabla v_{l_3}) \nonumber \\ 
   &\qquad\qquad\qquad\qquad\qquad\qquad\qquad\qquad\qquad+  (T^{jk} \partial_{x_j} u_{l_1} \partial_{x_k}u_{l_2})F_{T,\gamma_0}(u_{l_3}v_4+v_{l_3}u_{4})\Bigg\}, \nonumber
\end{align}
where $F_{T,\gamma_0}$ is a general term that depends smoothly only on $T^{jk}$ and $\gamma_0$, and their derivatives.
Following the same argument  as we did when proving $T^1+iT^2=0$ after \eqref{rewritesolu} in the $m=1$ case, 
we can prove the leading order term is $O(h^{-2})$ while other terms is $o(h^{-2})$. Hence we obtain
\begin{equation}\label{m2result1}
T^{11}+2iT^{12}-T^{22}=0.
\end{equation}
We similarly get 
\begin{align}\label{m2result2}
T^{11}-2iT^{12}-T^{22}=0
\end{align}
by choosing $u_1$, $u_2$, $u_3$ as in \eqref{ncsolu2} and $u_4$ as in $\eqref{csolu1}$.

Now we let the solutions $u_l$ be as in \eqref{ncsolu1} and \eqref{ncsolu2}:
\begin{align*}
   {u_1} & =\frac{1}{\sqrt{\gamma_0}}e^{(z+\frac{1}{2}z^2)/h}(1+r_1), \\
   {u_2} & =\frac{1}{\sqrt{\gamma_0}}e^{(-z+\frac{1}{2}z^2)/h}(1+r_2), \\
   {u_3} & =\frac{1}{\sqrt{\gamma_0}}e^{(-\bar{z}-\frac{1}{2}\bar{z}^2)/h}(1+r_3),\\
   {u_4} & =\frac{1}{\sqrt{\gamma_0}}e^{(\bar{z}-\frac{1}{2}\bar{z}^2)/h}(1+r_4).
\end{align*}
Again we  write $\partial_1 = \partial+\bar{\partial}$, $ \partial_2= i(\partial-\bar{\partial})$ and $\nabla u \cdot \nabla v = 2(\partial u \bar{\partial}v+ \bar{\partial} u \partial v)$
so that \eqref{eqT^ij_neq0} becomes
\begin{multline}\label{eqT11+T22}  
    \sum_{\left(l_1, l_2, l_3\right) \in \pi(3)}  \int \Big[T^{11}((\partial+\bar{\partial}) u_{l_1})( (\partial+\bar{\partial}) u_{l_2})+ T^{12} ((\partial+\bar{\partial})u_{l_1})( i(\partial-\bar \partial){u_{l_2}} )\\
    + T^{21} (i(\partial-\bar{\partial})u_{l_1}) ((\partial+\bar \partial){u_{l_2}}) -T^{22} ((\partial-\bar \partial)u_{l_1} )((\partial-\bar \partial)u_{l_2} ) \Big] (\partial u_{l_3}\bar{\partial}u_4+ \bar{\partial} u_{l_3} \partial u_4)  =0.
\end{multline}
Consider first $l_3=1$ in the above equation. 
Let us first focus on the terms in \eqref{eqT11+T22}  where $\partial$ hits  $u_1$ and $u_2$, and $\bar{\partial}$ hits $u_3$ and $u_4$. The resulting term is 
\begin{align}\label{quantity:largerterm}
    2(T^{11}+T^{22})(\partial u_1 \bar{\partial} u_3) 
(\partial u_{2} \bar{\partial}u_4).
\end{align}

Now if the derivatives in \eqref{quantity:largerterm} all hit the exponential factor of the solutions, we  get
\begin{multline*}
   2\int \frac{1}{h^4 \gamma_0^2} (T^{11}+T^{22}) e^{(z^2-\bar{z}^2)/h} (1+z)(-1+z)(1-\bar{z})(-1-\bar{z})\\
   \times (1+r_1)(1+r_2)(1+r_3)(1+r_4) = \frac{2}{h^3} ((T^{11}+T^{22})/\gamma_0)(0) + o(h^{-3}),
\end{multline*}
where we used estimate \eqref{ncestimate} and the same argument, that begun from \eqref{eq:first_integral}, involving the remainders in the case $m=1$.
Therefore, to show $T^{11}+T^{22}=0$,  it remains to show that in \eqref{eqT11+T22}  all the other terms are $o(h^{-3})$.

 We now consider in the above case $l_3=1$ that one derivative hits the $\gamma_0$ term. For example, consider the following term
$$
\int \frac{ \partial(1/\sqrt{\gamma_0})}{h^3\sqrt{\gamma_0}^3} T e^{(z^2-\bar{z}^2)/h} (-1+z)(1-\bar{z})(-1-\bar{z})(1+r_1)(1+r_2)(1+r_3)(1+r_4).
$$
In the expansion of the above of product, the term with no remainder term is of order $O(h^{-2})$ by stationary phase, and the terms with remainder term is also $O(h^{-2})$ by estimate \eqref{ncestimate}.
The other cases where one or more derivatives hit $\gamma_0$ term instead of the exponential term are similar. 

Now if one derivative hits the $(1+r)$ term, we consider
$$
\int \frac{1}{h^3 \gamma_0^2} T e^{(z^2-\bar{z}^2)/h} (-1+z)(1-\bar{z})(-1-\bar{z})(\partial r_1)(1+r_2)(1+r_3)(1+r_4).
$$
By estimate \eqref{ncestimate}, we  conclude that every term also in the above integral is $O(h^{-2})$. Similarly for all the other cases where derivatives hit $1+r$ term instead of the exponential term, we have the same conclusion.

Let us then consider the remaining terms in \eqref{eqT11+T22} where   $\partial$ does not hit both  $u_1$ and $u_2$ if  $\bar{\partial}$ hits both $u_3$ and $u_4$, and vice versa. For example, consider terms of form
\begin{align}\label{quantity:smallerterm}
    T\partial u_1 \partial u_3 \partial u_2 \bar{\partial} u_4.
\end{align}
Then the term where the derivative for $u_3$ hitting on the exponential will vanish. Therefore, we only need to consider 
$$
\int T e^{(z^2-\bar{z}^2)/h} \partial(1/\sqrt \gamma_0) \partial u_1  \partial u_2 \bar{\partial} u_4
$$
and 
$$
\int T e^{(z^2-\bar{z}^2)/h} (\partial r_3) \partial u_1  \partial u_2 \bar{\partial} u_4.
$$
By stationary phase, both integrals are $O(h^{-2})$.
The case $l_3=2$ is similar to $l_3=1$, thus we omit its proof.

Finally, we consider the case $l_3=3$. In this case, 
$(\partial u_{l_3} \bar{\partial}u_4+ \bar{\partial} u_{l_3} \partial u_4) 
$ becomes
$(\partial u_{3} \bar{\partial}u_4+ \bar{\partial} u_{3} \partial u_4) 
$, so we only need to consider terms of the form 
\begin{align}
    T\partial u_1 \partial u_3 \partial u_2 \bar{\partial} u_4 \qquad \text{or} \qquad T\partial u_1 \bar{\partial} u_3 \partial u_2 {\partial} u_4 .
\end{align}
Similar argument for \eqref{quantity:smallerterm} shows the above terms are also $O(h^{-2})$. Therefore, we have shown all the other terms except \eqref{quantity:largerterm} are $o(h^{-3})$, and
we get $T^{11}+T^{22}=0$ in $\Omega$. Combining with \eqref{m2result1} and \eqref{m2result2}, we conclude that $T^{11}=T^{12}=T^{22}=0$ in $\Omega$. This finishes the proof for $m=2$.

\vspace{12pt}

\subsection{The case $m\geq 3$:}
Let us recall the integral identity for general $m$: \eqref{prop1equ}
\begin{multline}\label{eq:mgeq3}
\sum_{\left(l_1, \ldots, l_{m+1}\right) \in \pi(m+1)} \sum_{j_1, \ldots, j_m=0}^2 \int T^{j_1 \cdots j_m}(x)\left(u_{l_1}, \nabla u_{l_1}\right)_{j_1} \ldots\left(u_{l_m}, \nabla u_{l_m}\right)_{j_m} \\
\times \nabla u_{l_{m+1}} \cdot \nabla u_{m+2}\s =0.
\end{multline}
We first prove $T^{\s0\s j_1\s \cdots\s j_{m-1}}=0$ where $j_k\in \{0,1,2\}$ and $k=1,\ldots,m-1$.
Firstly, let $m$  of the functions $u_l$, $l=1,\ldots,m+1$, in the identity \eqref{eq:mgeq3} to be the constants functions $1$.
%
This yields 
\begin{equation}\label{eq:red_to_m1}
    \int  T^{\s0\s\cdots\s0}(x) \nabla u \cdot \nabla v  = 0, 
\end{equation}
which is of the form \eqref{s3e1} we had in the case $m=1$. Thus  we obtain  $T^{\s0\s\cdots\s 0}=0$ in $\Omega$. 

Next, we let $m-1$ of solutions $u_l$  to be the constant functions $1$. This yields
\begin{equation}
    \sum_{\left(l_1, l_2\right) \in \pi(2)} \sum_{j=1}^2 \int T^{0\ldots0j}(x) \partial_{x_j} u_{l_1} \nabla u_{l_2} \cdot \nabla u_3  = 0,
\end{equation}
which is of the form \eqref{eq3} we also had in the case $m=1$. Thus we have  $T^{\s0\s\cdots\s 0\s j}=0$ in $\Omega$ for $j=1,2$. Continuing in similar fashion, we let $m-2$ of the solutions $u_i$ to be the constant functions $1$. This yields an integral identity similar to \eqref{eqT^ij_neq0} we had in the case $m=2$. The same argument used in that case proves
$T^{\s0\s\cdots\s 0\s j_1\s j_2}=0$ in $\Omega$ for $j_1,j_2\in \{1,2\}$.
Proceeding in this manner, by induction we obtain 
\begin{equation}\label{eq:induc_on_T}
T^{\s0\s j_1\s \cdots\s j_{m-1}}=0, \quad j_1,\ldots, j_{m-1} \in \{0,1,2\}
\end{equation}
in  $\Omega$.

It remains to prove $T^{j_1\cdots j_m}=0$ where all the indices $j_k$, $k=1,\ldots, m$, are nonzero. By \eqref{eq:induc_on_T}, the integral identity \eqref{eq:mgeq3} is reduced to
\begin{align}\label{eq10}
    \sum_{\left(l_1, \cdots, l_{m+1}\right) \in \pi(m+1)} \sum_{j_1,\ldots j_m=1}^2 \int T^{j_1 \ldots j_m}(x) \partial_{x_{j_1}} u_{l_1} \partial_{x_{j_2}} u_{l_2} \ldots \partial_{x_{j_m}} u_{l_m} \nabla u_{l_{m+1}} \cdot \nabla u_{m+2} =0
\end{align} for all $u$ solving \eqref{cond_equ}.
Since $T^{j_1\cdots j_m}$ is symmetric in exchange of any of its two indices, it has $(d+m-1)!/(m!(d-1)!)$ independent components, where $d=2$. Thus we have $m+1$ unknown entries $T^{j_1\cdots j_m}$ in \eqref{eq10}. To recover these entries we will find $m+1$ linearly independent  equations for the entries.
Firstly, we choose CGO solutions such that exactly one of them has an antiholomorphic phase:
\begin{align*}
   {u_1}=u_2=\cdots=u_m & =\frac{1}{\sqrt{\gamma_0}}e^{(z+\frac{1}{m+1}z^2)/h}(1+r_1), \\
   {u_{m+1}} & =\frac{1}{\sqrt{\gamma_0}}e^{(-mz+\frac{1}{m+1}z^2)/h}(1+r_2), \\
   {u_{m+2}} & =\frac{1}{\sqrt{\gamma_0}}e^{-\bar{z}^{2}/h}(1+r_3).
\end{align*}
Following how we proved $T^1+iT^2=0$ in $m=1$ case, we can show that the principal order term in \eqref{eq10} is $O(h^{-m})$, and the integrals involving correction terms $r_1, r_2, r_3$ are $o(h^{-m})$. Stationary phase shows that the principal order term of \eqref{eq10} gives the linear equation 
\begin{align}\label{eq:coefficientline1}
   \binom{m}{0} T^{1\ldots1}+i\binom{m}{1}T^{1\ldots12}+i^2 \binom{m}{2}T^{1\ldots122}+\ldots+i^{m}\binom{m}{m}T^{2\ldots2}&=0,
\end{align}
after dividing by a nonzero constant. To see how \eqref{eq:coefficientline1}\ is obtained, we first note that the coefficients come from expanding  the integral
\begin{equation}\label{eq:integral_id_for_gen_m}
   \sum_{j_1,\ldots j_m=1}^2 \int \bar{\partial} T^{j_1 \cdots j_m}(x) (\partial \partial_{x_{j_1}} u_{l_1}) (\partial_{x_{j_2}} u_{l_2}) \ldots (\partial_{x_{j_m}} u_{l_m}) (u_{l_{m+1}} u_{m+2})
\end{equation}
 by stationary phase. 
Here the principal order term results from the solutions $\partial$ hitting $u_1, u_2, \ldots, u_{m+1}$ and $\bar{\partial}$ hitting the solution $u_{m+2}$, which have holomorphic and antiholomorphic phases respectively. This is similar to what we had in the proof for the case $m=1$. We also used 
\begin{equation}
\partial \partial_{1} = \partial^2+\partial \bar{\partial}, \quad
\partial \partial_{2} = i(\partial^2-\partial \bar{\partial}).
\end{equation}
to compute the exact coefficients.


Similarly, if we choose CGO solutions such that exactly one of them has holomorphic phase, we get the following linear equation up to a scalar multiple: 
\begin{align}\label{eq:coefficientline2}
   \binom{m}{0} T^{\s1\s\cdots\s1}+(-i)\binom{m}{1}T^{\s1\s\cdots\s1\s2}+(-i)^{2} \binom{m}{2}T^{\s1\s\cdots\s1\s2\s2}+\cdots+(-i)^{m}\binom{m}{m}T^{\s2\s\cdots\s2}&=0.
\end{align}

Next we choose CGO solutions so that more than one of the solutions have holomorphic phases and also that more than one solution have antiholomorphic phases. Therefore, we can choose every solution to be of form \eqref{ncsolu1} or \eqref{ncsolu2} whose phase functions have no critical points. For general $m$, we choose CGOs to that all their phases add up to $z^2 - \bar{z}^2$. The explicit formula 
 for general $m$ is complicated to write down.
 Therefore, we only consider the case $m=3$ as an example. In this case, we choose two solutions with holomorphic phases and three with antiholomorphic phases: 
\begin{align*}
   {u_1} & =\frac{1}{\sqrt{\gamma_0}}e^{(z+\frac{1}{2}z^2)/h}(1+r_1), \\
   {u_2} & =\frac{1}{\sqrt{\gamma_0}}e^{(-z+\frac{1}{2}z^2)/h}(1+r_2), \\
   {u_3} & =\frac{1}{\sqrt{\gamma_0}}e^{(-\bar{z}-\frac{1}{3}\bar{z}^2)/h}(1+\tilde{r_3}),\\
   {u_4} & =\frac{1}{\sqrt{\gamma_0}}e^{(-\bar{z}-\frac{1}{3}\bar{z}^2)/h}(1+\tilde{r_4}),\\
   {u_5} & =\frac{1}{\sqrt{\gamma_0}}e^{(2\bar{z}-\frac{1}{3}\bar{z}^2)/h}(1+\tilde{r_5}).
\end{align*}
Note that here all the solutions have phases without critical points. Consequently, their correction terms satisfy the better estimates \eqref{ncestimate}, which simplifies the asymptotic analysis.

By arguing similarly as we did after \eqref{eqT11+T22}, we obtain
\begin{align}\label{1-i1-i}
    T^{111}-iT^{112}+T^{122}-iT^{222}&=0.
\end{align}
Note that the coefficients in \eqref{1-i1-i} agree with those in the expansion of the polynomial $(a+ib)(a-ib)^2$ of variables $a$ and $b$. This is true in general: 
Choose $t$, $1\leq t \leq m+1$, solutions to have holomorphic phases and $m+2-t$ solutions to have antiholomorphic phases in \eqref{eq:mgeq3}. 
Then, by stationary phase, we may compute the coefficient of $T^{\s1\s\cdots\s1\s2\s\cdots2}$, where the number of indices with index $1$ is $s$  and the number of indices with index $2$ is $(m-s)$. 
The coefficient will agree with the coefficient of $a^s b^{m-s}$ in the expansion of $(a+ib)^{t-1} (a-ib)^{m+2-t-1}$. We explain next why the above holds.

The reason why the above holds is the following: 
The principal order term of the integral \eqref{eq:integral_id_for_gen_m}  for the chosen solutions corresponds to  $\partial$ hitting solutions with holomorphic phases and $\bar{\partial}$ hitting solutions  with antiholomorphic phases. (See the part of the proof after \eqref{eqT11+T22}.) Then, since
$$
\partial_1 = \partial+\bar{\partial}, \quad \partial_2 = i(\partial-\bar{\partial}),
$$ 
we know that $\partial_1$ acting on a holomorphic phase gives the coefficient  $1$, $\partial_2$ acting on holomorphic phase gives a coefficient of $i$, $\partial_1$ acting on antiholomorphic phase gives a coefficient of $1$ and $\partial_2$ acting on antiholomorphic phase gives a coefficient of $-i$.

To compute the coefficient of  $T^{\s 1\s \cdots\s1\s2\s\cdots\s2}$, where the number of indices $1$ is $s$, and the number of indices $2$ is $(m-s)$, we note the following.  In the integral \eqref{eq10}, the coefficient $T^{\s 1\s \cdots\s1\s2\s\cdots\s2}$ appears together with $s$ instances of $\partial_1$ and $m-s$ instances of $\partial_2$. In the principal order term, if $u$ is a solution with holomorphic phase, we may only consider the terms where we have $\partial u$. So for the term $\partial_1 u = (\partial + \bar{\partial})u$ following $T^{\s 1\s \cdots\s1\s2\s\cdots\s2}$ in \eqref{eq10}, we may consider only $\partial u$, while in the term $\partial_2 u = i(\partial - \bar{\partial})u$, we may consider only $i\partial u$. 

We have similarly for $u$ with antiholomorphic phase: For $\partial_1 u = (\partial + \bar{\partial})u$, we consider only  $\bar{\partial} u$, while in the term $\partial_2 u = i(\partial - \bar{\partial})u$, we consider only  $-i\bar{\partial} u$. This implies each of $\partial_1$ corresponds to a factor of $1$ for both holomorphic phase and antiholomorphic phase, and each $\partial_2$ corresponds to a factor of $i$ for holomorphic phase, and a factor of $-i$ for antiholomorphic phase.   From the proof after $\eqref{eqT11+T22}$, we see there are $t-1$ solutions having holomorphic phases and $m+2-t-1$ solutions having antiholomorphic phases among $u_{l_1},\ldots, u_{l_m}$. Therefore, the coefficient that we are considering should be equal to that of $a^s b^{m-s}$ in the polynomial $(a+ib)^{t-1} (a-ib)^{m+2-t-1}$.

Note that in the case of $m$, we have $m+1$ many  choices of solutions. Each choice gives a linear equation of the form \eqref{eq:coefficientline1}. Finally, we show that the $m+1$ linear equations we have now obtained for the coefficients $T^{j_1j_2j_3}$ are linearly independent. This implies that the coefficients are uniquely determined. Let us inspect the linear system we obtain for $m=3$. This is 
\begin{align*}
    T^{111}+3iT^{112}-3T^{122}-iT^{222}&=0\\
    - T^{111}-iT^{112}-T^{122}-iT^{222}&=0\\
    T^{111}-iT^{112}+T^{122}-iT^{222}&=0\\
    -T^{111}+3iT^{112}+3T^{122}-iT^{222}&=0.
\end{align*}
As we have shown above, the coefficients in each row agree with those of $(a+ib)^3$, $(a+ib)^2(a-ib)$, $(a+ib)(a-ib)^2$ and $(a-ib)^3$ respectively. Since these polynomials are linearly independent, so is the coefficient matrix of the above linear system. The proof for general  $m$ is similar.
\\ \vspace{6pt}

\bibliographystyle{alpha}
\bibliography{2100_v2}


\end{document}